\documentclass[12pt,oneside]{amsart}
\usepackage{tkz-euclide}
\usepackage{tkz-euclide}
\usetikzlibrary{positioning}
\usepackage{geometry}                
\usepackage{graphicx}
\usepackage{amssymb}
\usepackage{epstopdf}

\usepackage{amsmath,amsthm,amscd,amssymb, enumerate}
\usepackage{latexsym}
\usepackage[colorlinks,citecolor=red,pagebackref,hypertexnames=false]{hyperref}
\numberwithin{equation}{section}
\theoremstyle{plain}
\newtheorem{theorem}{Theorem}[section]

\newtheorem{corollary}[theorem]{Corollary}
\newtheorem{proposition}[theorem]{Proposition}

\theoremstyle{definition}
\newtheorem{definition}[theorem]{Definition}

\theoremstyle{remark}
\newtheorem{remark}[theorem]{Remark}

\newtheorem{case[theorem]}{Case}

\def\bc{\begin{corollary}}
\def\ec{\end{corollary}}
\def\be{\begin{equation}}
\def\ee{\end{equation}}
\def\bast{\begin{eqnarray*} }
\def\east{\end{eqnarray*} }
\def\bea{\begin{eqnarray}}
\def\eea{\end{eqnarray}}

\def\bmat{\begin{matrix}}
\def\emat{\end{matrix}}

\def\ms{\medskip}
\def\bs{\bigskip}

\def\R{\mathbb R}
\def\hd{{\dim_{\mathcal H}}}

\def\R{\Bbb R}

\def\vt{{\bf{t}}\,}

\def\dt{d_{\text{tot}}}
\def\dl{d_L^\sigma}
\def\dr{d_R^\sigma}
\def\meff{m_{\text{eff}}}
\def\pk{\mathcal P_k}

\def\hd{\dim_{\mathcal H}}

\def\({\left(}
\def\){\right)}
\def\[{\left[}
\def\]{\right]}
\def\<{\left\langle}
\def\>{\right\rangle}

\def\zvt{Z_{\vt}}
\def\zvts{\zvt^{\sigma}}
\def\rvt{\mathcal R_{\vt}}

\def\rvts{\rvt^\sigma}
\def\rvtsl{\rvt^{\sigma^l}}
\def\rp{\R\mathbb P^2}

\def\sd{\mathbb S^{d-1}}

\def\om{\omega}
\def\tom{\tilde{\om}}
\def\ttom{\om'}
\def\ti{\tilde{i}_{\tom}}

\def\det{\hbox{det}}
\def\dim{\hbox{dim}}

\def\inter{\hbox{\emph{Int}}\,}
\def\rank{\hbox{rank }}

\def\wtnzt{\widetilde{N^*\zvt}}

\def\tu{\tilde{u}}
\def\tw{\tilde{w}}
\def\ty{\tilde{y}}
\def\tz{\tilde{z}}

\def\atxt{}
\def\aatxt{}

\thanks{Keywords: Erd{o}s-Falconer configuration problem, Radon transform, Fourier integral operator. 
MSC:  28A75, 28A80, 58J40 (primary), 52C10 (secondary).}
\thanks{
AG is supported in part by  US National Science Foundation DMS-1906186 and -2204943, AI by NSF HDR TRIPODS-1934962 and DMS-2154232, 
and  KT   by  Simons Foundation Grant 523555.
{The authors thank the referee for helpful suggestions and questions which improved the paper.}}

\begin{document}

\author{Allan Greenleaf}
\address{Department of Mathematics, University of Rochester, Rochester, NY 14627}
\email{allan@math.rochester.edu}
\author{Alex Iosevich} 
\address{Department of Mathematics, University of Rochester, Rochester, NY 14627}
\email{iosevich@math.rochester.edu}
\author{Krystal Taylor}
\address{Department of Mathematics, Ohio State University, Columbus, Ohio 43210}
\email{taylor.2952@osu.edu}

\title[Configuration sets via microlocal partition optimization]
{Nonempty interior of configuration sets \\ via microlocal partition optimization}

\date{Revision of  February 16, 2024. 
}    

\maketitle

\begin{abstract} We prove new results of Mattila-Sj\"olin  type,
giving lower bounds on  Hausdorff dimensions
 of  thin sets  $E\subset \R^d$ ensuring that various $k$-point configuration sets,  
generated by elements of $E$,  have nonempty  interior. 
The dimensional thresholds in our previous work \cite{GIT20}
were dictated by associating to a configuration function
a family of generalized Radon transforms, and then 
optimizing $L^2$-Sobolev  estimates for them over all nontrivial bipartite partitions of the $k$ points.
In the current work, we extend this by allowing the optimization to be  
done locally over the configuration's incidence relation, 
or even microlocally over   the conormal bundle of the incidence relation. 
\aatxt{We use this approach  to 
prove Mattila-Sj\"olin  type results for 
(i) areas of subtriangles 
determined  by quadrilaterals and pentagons
 in a set $E\subset\R^2$;
(ii) pairs of ratios of distances of 4-tuples in $\R^d$; and
 (iii)  similarity classes of triangles in $\R^d$,
 as well as to  
 (iv) give
 a short  proof  of Palsson and Romero Acosta's result  
on congruence classes of triangles in $\R^d$.}
\end{abstract}  

\maketitle


\section{Introduction} \label{sec intro}

\aatxt{This is the third in a series of papers using Fourier integral operator techniques to obtain Mattila-Sj\"{o}lin type results, by 
which we mean results showing that certain types of configuration sets have nonempty interior when the underlying sets have 
sufficiently high Hausdorff dimension.
In  \cite{GIT19,GIT20}, we showed how to obtain such results
for a wide range of configurations, using estimates for generalized Radon transforms,
based on analysis of them as linear or multilinear Fourier integral operators.
The current paper extends these methods and gives several applications. }
\ms

A classical result of Steinhaus \cite{Steinhaus} states that if $E\subset \R^d, d\ge 1,$ has positive Lebesgue measure, 
then the difference set $E-E\subset \R^d$ contains a neighborhood of the origin.
 $E-E$ can be {identified with}  the set of two-point configurations, $x-y$, of points of $E$ modulo the translation group.
In the context of the Falconer distance set problem,  
a  theorem of Mattila and Sj\"olin \cite{MS99} states that if $E\subset\R^d,\, d\ge 2$, is compact,
then the distance set of $E$, 
 $\Delta(E)=:\{\, |x-y|: x,\, y\in E\}\subset\R$,
 contains an open interval,  i.e., has nonempty interior,
if  the Hausdorff dimension $\hd(E)>\frac{d+1}2$. 
This provided a strengthening of Falconer's original result \cite{Falc86},
from $\Delta(E)$ merely having positive Lebesgue measure to having nonempty interior,  for the same range of $\hd(E)$.
 This was generalized   to distance sets  with respect to  norms on $\R^d$ 
 having positive curvature unit spheres in  Iosevich, Mourgoglou and Taylor \cite{IMT11}.
 \ms

Mattila-Sj\"olin type results, establishing nonempty interior for sets of 
configurations in a set $E$
 only satisfying a lower bound on $\hd(E)$,
or results that can be interpreted as such,
have been obtained by various authors.
These include \cite{IMT11,GIP14,BIT16,CLP16,IT19} and, more recently, \cite{KPS21,PRA21,McDT21};
see also \cite{CKPS21} for a finite field analogue, as well as \cite{McDT21, McDT22, SiT20} for analogues in which Hausdorff 
dimension is replaced by an alternative notion of size, mainly Newhouse thickness.  
\ms

More general  Mattila-Sj\"olin style theorems were studied by the current authors,  for $2$-point configurations in \cite{GIT19} 
and $k$-point configurations in \cite{GIT20}.
In those, as in the present work, the configurations considered are  $\Phi$\emph{-configurations}, 
as defined   by Grafakos, Palsson and the first two authors  \cite{GGIP12},
which can be vector-valued, nontranslation-invariant
and possibly  {asymmetric}, i.e.,  among points in sets $E_1,\dots, E_k$ lying in 
different spaces, e.g., points and  circles in $\R^2$. 
The approach taken was to study the $L^2$-Sobolev mapping properties
of an associated family of generalized Radon transforms, linear in \cite{GIT19} or multilinear in \cite{GIT20}.
The main step in showing that the set $\Delta_\Phi(E_1,\dots,E_k)$ has nonempty interior is analysis of the configuration  measure  $\nu(\vt)$ (defined below); we show that this measure is absolutely  continuous and that its density with respect to Lebesgue measure, $d\vt$, is a continuous function of the configuration parameter $\vt\in\R^p$ (or other space). 
This was done in \cite{GIT20} by representing $\nu(\vt)$ as the pairing of the tensor product 
of  Frostman measures $\mu_i$ on some of the $E_i$ 
with the value of a generalized Radon  transform, $\rvt$,
acting on the tensor product of Frostman measures on the complementary collection 
of $\mu_j$-s.
\aatxt{Each such partition of the $k$ variables into two  groups gives
a threshold for $\sum_{j=1}^k \hd(E_j)$ ensuring 
$\inter\left(\Delta_\Phi(E_1,\dots,E_k)\right)\ne\emptyset$; the threshold can potentially be lowered by optimizing over all 
such partitions.}
 We  refer to that approach  as \emph{partition optimization}; for a precise 
statement, see Theorem \ref{thm kpoint old} below, 
which is \cite[Thm. 5.2]{GIT20}. \aatxt{(See \cite{GaGrPaPs24} for a subsequent application of partition optimization.)}
\ms

The purpose of the current paper is to show that  extensions of that approach,
{performing the partition optimization 
\emph{locally} or even  \emph{microlocally},}
allow one to obtain such nonempty interior results for an even wider range of $k$-point 
configurations,  which fail to satisfy the hypotheses of Theorem \ref{thm kpoint old}. 
We do this by considering open covers
of the  $k$-fold incidence relation defining the configuration of interest,
or more generally allowing
microlocal  covers of the conormal bundle of the incidence relation
by open, conic sets.
\aatxt{(This microlocal refinement of the local method is not possible for incidence relations of codimension one, since above 
each point the incidence relation is just a line, and an open conic cover of the conormal bundle is equivalent to an open cover of 
the incidence relation.)}
On each  of these  sets, Theorem \ref{thm kpoint old} is applicable,
but with the Hausdorff dimensional threshold  \aatxt{\it possibly optimized by  different} partitions of the 
$k$ variables as one ranges over the elements of the cover.
Taking the maximum of the thresholds needed,
either locally near each point of the incidence relation
or microlocally near each point of its  conormal bundle, 
and then optimizing over all covers,  
yields a threshold which is always less than or equal to 
that provided by Theorem \ref{thm kpoint old};
see Theorems \ref{thm kpoint new local} and \ref{thm kpoint new} for the statements of the local and microlocal versions. 
See Sec. \ref{sec framework} for the background
material from \cite{GIT20} and the precise statements and the proofs of the theorems.
\smallskip

We now state some results  which can be obtained using this new approach,
restricting the discussion to various three-, four- and five-point configurations in $\R^d$.
\ms

{\it Areas of triangles generated by vertices of quadrilaterals and pentagons in $\R^2$:} In \cite[Thm. 1.1]{GIT20} we 
showed that if $E\subset\R^2$ with $\hd(E)>5/3$, then the set of areas of 
triangles with vertices in $E$ has nonempty interior in $\R$. 
\aatxt{For $n$-tuples of vertices in $E$, with $n\ge 4$, one can also consider vector-valued configurations consisting of
the areas of {\it some} of the triangles they generate. }
In \cite{GIT20} we established that the collection of ordered pairs of 
areas of two of the triangles  generated by a quadrilateral $xyzw$ 
with vertices in $E$, say $\left(|xyz|,|xzw|\right)$, has nonempty interior in $\R^2$ 
if $\hd(E)>7/4$. Note that there are  limits 
on how far such results can be pushed: since 
$|xyw| + |yzw|=|xyz| +|xzw|$, the configuration set 
of all four of these areas would 
lie in a hyperplane in $\R^4$ and thus would have empty interior. However, using 
microlocal partition optimization, we are able to obtain  \linebreak
(i) a threshold 
improving upon that in \cite[Thm. 1.6]{GIT20}; 
(ii) a result for triples of areas of triangles generated by a quadrilateral; 
and (iii) a result for triples of the areas of a fan of triangles generated by a pentagon.
\ms

\begin{theorem}\label{thm triangles}
If $E\subset\R^2$ is compact, then
\ms

(i) if $\hd(E)>3/2$, then $\inter\left\{(|xyz|,|xzw|)\in\R^2: x,y,z,w\in E\right\}\ne\emptyset$; 
\ms

(ii) if $\hd(E)>7/4$, then $\inter\left\{(|xyz|,|xzw|,|xyw|)\in\R^3: x,y,z,w\in E\right\}\ne\emptyset$; and
\ms

(iii) if $\hd(E)>9/5$, then $\inter\left\{ (|xyz|,|xzw|,|xwu|)\in\R^3: x,y,z,w,u\in E \right\}\ne\nolinebreak\emptyset$.
\end{theorem}
\bs

{For the proof of Theorem \ref{thm triangles}, see Sec. \ref{sec triangles}. Also see A. McDonald \cite{McD21} for related results of 
Falconer  type (i.e., positive Lebesgue measure), 
formulated in terms of areas of parallelograms generated by pairs of points  in $E$ rather than areas of 
triangles generated by triples. }
\ms

{{\it Pairs of ratios of  distances:} 
Another result that can be obtained using
microlocal partition optimization concerns 
ratios of distances of 4-tuples in a set.
}
\ms

\begin{theorem}\label{thm ratios}
{If $E\subset \R^d,\, d\ge 2$, is  compact  and $\hd(E)>(3d+1)/4$, then 
$$\inter\left\{\left(\frac{|x-y|}{|z-w|},\, \frac{|x-w|}{|z-y|}\right)\in\R^2: 
x,y,z,w\in E, z\ne w,\, z\ne y \right\}\ne\emptyset.$$
}
\end{theorem}

\aatxt{In contrast, a  modification of this configuration, \eqref{eqn zywy}, has nonempty interior} for a larger range of dimensions
and has a more elementary proof.
A discussion of these and
other results concerning configuration sets defined by ratios of distances, can be found in Sec. \ref{sec ratios}. 
\bs

{\it Congruences classes of triangles in $\R^d,\, d\ge 4$:}
One motivation for developing the microlocal extension of  the original partition optimization technique  
was from trying to understand how a recent result of
Palsson and Romero-Acosta \cite{PRA21} related to the FIO framework of \cite{GIT20}.
They  proved the following:

\begin{theorem}\label{thm pra} {\bf \cite{PRA21}.}
If $E\subset\R^d,\, d\ge 4,$ is compact with $\hd(E)>(2d+3)/3$,
then  the set of congruence classes of triangles with vertices in $E$,
\be\label{eqn pra}
\left\{\, \(\,|x-y|,\, |x-z|,\, |y-z|\, \):\, x,y,z\in E\, \right\},
\ee
 has nonempty interior in $\R^3$. 
\end{theorem}

{In Sec. \ref{sec PRA}, we present a much shorter proof of this.
Both the original partition optimization 
method from \cite{GIT20}  (Thm. \ref{thm kpoint old} below)  and the local version 
(Thm. \ref{thm kpoint new local} below) } fail to prove Thm. \ref{thm pra},  
\aatxt{because the optimal partition of the three variables varies with the normal direction to the (codimension 3) incidence relation, even above a single point.}
However, it can be proved using  \aatxt{microlocal  partition optimization}, Thm. \ref{thm kpoint new}.  
\ms

\atxt{{\it Similarity classes of triangles in $\R^d,\, d\ge 3$.} We conclude with a result that has 
similarities to both Thm. \ref{thm ratios} and Thm. \ref{thm pra}. 
Using Thm. \ref{thm kpoint new local},  in Sec. \ref{sec similarity} we prove}

\atxt{\begin{theorem}\label{thm similarity} 
If $E\subset\R^d,\, d\ge 3$, is compact and $\hd(E)>(2d+2)/3$, then the set of similarity classes of triangles with vertices
in $E$, 
\be\label{eqn sim}
\Big\{ \big[\,  |x-y|\, :\, |x-z|\, :\, |y-z|\, \big]\in \rp \hbox{ s.t. } x,y,z\in E\hbox{ are distinct }\Big\}
\ee
has nonempty interior in $\rp$. 

\aatxt{(Here, $[\, A:\, B\, :\, C\, ]$ are standard projective coordinates on $\rp$.)}
\end{theorem}
}
\bs

\atxt{Before proving the theorems, in Sec. \ref{sec framework} we  recall  the framework of $\Phi$-configurations 
and the method of partition optimization   from \cite{GIT19,GIT20}.}

 \section{$k$-point 
 $\Phi$-configuration sets}\label{sec framework}
 
 In order to state microlocal partition optimization, \aatxt{and for the sake of readability,} 
 we  recall from  \cite{GIT19,GIT20} the framework 
 for studying the $\Phi$-configuration sets of
 \cite{GGIP12} via FIO methods \aatxt{ and the original (global) version of partition optimization.} 
Suppose that  $X^i,\, 1\le i\le k,$ and $T$, are  
 smooth manifolds of dimensions $d_i$ and $p$, 
 resp. We sometimes denote  $X^1\times\cdots\times \nolinebreak X^k$ by $X$, and set $\dt:=\dim(X)=\sum_{i=1}^k d_i$.
 
 \begin{definition}\label{def k Phi config} 
  Let $\Phi\in C^\infty(X,T)$. 
Suppose that   $E_i\subset X^i,\, 1\le i\le k$, are compact sets.
Then the {\it $k$-configuration set of the $E_i$} defined by $\Phi$ is
\be\label{def config general}
\Delta_\Phi\left(E_1,E_2,\dots,E_k\right):=\left\{\Phi\left(x^1,\dots,x^k\right): x^i\in E_i,\, 1\le i\le k\right\}\subset T.
\ee
If $E=E_1=\cdots=E_k=E$, then we just write $\Delta_\Phi(E)$.
\end{definition}

\medskip

We want to find  conditions on  the $\hd(E_i)$ ensuring that  $\Delta_\Phi\left(E_1,E_2,\dots,E_k\right)$
has nonempty interior.
To this end, now suppose that $\Phi:X\to T$ is a submersion, so that for each $\vt\in T$, $\zvt:=\Phi^{-1}(\vt)$ is a smooth, 
codimension $p$ submanifold of $X$, 
and these vary smoothly with $\vt$. 
For each $\vt$, the measure
 \be\label{def delta k}
\lambda_\vt:= \delta\(\Phi\(x^1,\dots,x^k\)-\vt\)
 \ee
 is a smooth density on $\zvt$; i.e., a smooth multiple of surface measure. In local coordinates {$\vt=(t_1,\dots,t_p)$} on $T$, $\lambda_\vt$  can be represented as an oscillatory integral
of the form
 \be\label{eqn multi ker}
\lambda_\vt= \int_{\R^p} e^{i\[\sum_{l=1}^p\(\Phi_l(x^1,\dots,\, x^k)-t_l\)\tau_l\]}\, a(\vt)1(\tau)\, d\tau,
 \ee
 where the $a(\cdot)$ belongs to a partition of unity on $T$,
 { and $1(\tau)$ comes from the Fourier transform of the delta distribution in $\R^p$.}
 Thus, $\lambda_\vt$ is a Fourier integral distribution on $X$; 
 in H\"ormander's notation \cite{Hor71,Hor85,GuSt},
 \be\label{eqn fid}\nonumber
 \lambda_\vt\in I^{(2p-\dt)/4}(X;N^*\zvt),
 \ee
  where $N^*\zvt\subset T^*X\setminus 0$ is the conormal bundle of $\zvt$ and
  the value of the  order  follows
from the amplitude having order zero and the numbers of phase variables and spatial variables  being $p$ and $\dt$, resp., 
so that the order is $m:=0+p/2-\dt/4$.
\ms

We separate the variables $x^1,\,\dots, x^k$ into groups on the left and right,
  associating to $\Phi$ a collection of families of generalized Radon transforms  
indexed by the nontrivial partitions of $\{1,\dots,k\}$, 
with each family then 
depending on  the parameter $\vt\in T$.
Write such a partition   as $\sigma=(\sigma_L\, |\, \sigma_R)$,
with $|\sigma_L|,\, |\sigma_R|>0,\, |\sigma_L|+|\sigma_R|=k$, 
and let $\pk$ denote the set of all $2^k-2$ such partitions.
We use $i$ and $j$ to refer to elements of $\sigma_L$ and $\sigma_R$, resp.
Define $\dl=\sum_{i\in \sigma_L} d^i$ and $\dr=\sum_{j\in \sigma_R} d^j$, so that $\dl+\dr=\dt$.
\ms

For each $\sigma\in\pk$, $\sigma_L=\left\{i_1,\dots,i_{|\sigma_L|}\right\}$ 
and $\sigma_R=\left\{j_1,\dots, j_{|{\sigma_R}|}\right\}$,
where without loss of generality we may assume that $i_1<\cdots < i_{|\sigma_L|}$ and $j_1<\cdots< j_{|\sigma_R|}$.
With a slight abuse of notation
we still denote  the coordinate-partitioned version of $x$ as $x$,
$$x=\left(x_L;x_{R}\right):=\(x^{i_1},\dots,x^{i_{|\sigma_L|}};\,  x^{j_1},\dots,x^{j_{|\sigma_R|}}\).$$

Write the  corresponding reordered  Cartesian product as
$$X_L \times X_{R} := \(X^{i_1} \times \cdots \times X^{i_{|\sigma_L|}}\) \times
\(X^{j_1} \times \cdots \times X^{j_{|\sigma_R|}}\);$$
again by abuse of notation, we sometimes still refer to this as $X$.
The dimensions of the two factors are $dim(X_L)=\dl$ and $dim(X_R)=\dr$, resp.
The choice of $\sigma$  also  defines a coordinate-partitioned version of each $\zvt$, 
\be\label{def zsigma}
\zvts:=\left\{ \(x_L; x_{R}\):\, \Phi\(x\)=\vt  \right\}\subset X_L \times X_{R},
\ee
with spatial projections to the left and right,  $\pi_{X_L}:\zvts\to X_L$ and $\pi_{X_R}:\zvts\to X_{R}$.
The  integral geometric double fibration condition for $\zvts$ is the requirement that
\be\label{def DF general}
 (DF)_{\sigma}\qquad \pi_{L}:\zvts\to X_L\hbox{ and }\pi_{R}:\zvts\to X_{R}\hbox{ are 
 submersions.}
\ee 
(See \cite{Helgason,GuSt,GuSt79}.) Note that, for a given $\sigma$,
a necessary (but not sufficient) condition for   $(DF)_\sigma$ to hold  is $p\le \dl\wedge \dr:= \min\left(\dl, \dr\right)$.
\medskip

If $(DF)_\sigma$  holds, then the generalized Radon transform $\rvts$,
defined weakly by
$$\rvts f(x_L)=\int_{\left\{x_{R}:\, \Phi\left(x_L,x_{R}\right)=\vt\right\}} f(x_{R}),$$
where the integral is with respect to the surface measure induced by $\lambda_\vt$ on the codimension $p$ submanifold
$\left\{x_{R}:\, \Phi\left(x_L,x_{R}\right)=\vt\right\}=
\left\{x_{R}: \left(x_L,x_{R}\right)\in \zvts\right\} \subset X_{R}$,
which extends from mapping $\mathcal D\left(X_{R}\right)\to \mathcal E\left(X_L\right)$ to
$$\rvts:\mathcal E'(X_{R})\to \mathcal D'(X_L).$$
{Here, $ \mathcal E,\, \mathcal D$ are the standard spaces of $C^\infty$ functions and those of compact support, resp., and $\mathcal E',\, \mathcal D'$ their dual spaces of distributions.} 
Furthermore,  
\be\label{def ctsigma}
C_{\vt}^\sigma:=\(N^*\zvts\)'=\left\{\(x_L,\xi_L;\, x_{R},\xi_{R}\):
\, \(x_L,x_{R}\)\in\zvts,\, \(\xi_L,-\xi_{R}\)\perp T\zvts \right\}
\ee
is contained in $\(T^*X_L\setminus 0\)\times \(T^*X_R\setminus 0\)$. 
Thus, $\rvts$ is an FIO,  $\rvts\in I^m\(X_L, X_R; C_{\vt}^\sigma\)$,
where the order $m$ is determined as in \eqref{eqn fid}  by  $m=0+p/2-\dt/4$ \cite{Hor71,Hor85}. 
Given the possible difference in the dimensions of $X_L$ and 
$X_R$,  due to the clean intersection calculus it is useful to express $m$  as
$$m=\meff^\sigma-\frac14\left|\dl-\dr\right|,$$
where the {\it effective} order of $\rvts$  is defined to be
\be\label{eqn msigma}
\meff^\sigma:=\(2p-\dt+\left|\dl-\dr\right|\)/4=\left(p-\left(\dl\wedge\dr\right)\right)/2.
\ee
By standard estimates for FIO \cite{Hor71,Hor85}, if $C_{\vt}^\sigma$  is a nondegenerate canonical relation, i.e., the  cotangent space projections
$\pi_L: C_{\vt}^\sigma \to T^*X_L$
and 
$\pi_R:C_{\vt}^\sigma \to T^*X_R$
have differentials of  maximal rank, then 
$$\rvts:L^2_r\left(X_R\right)\to L^2_{r-{\meff}^\sigma}\left(X_L\right).$$
More generally, if $\pi_L$ (and thus $\pi_R$) drops rank by $\le q$, then there is a loss of $\le q/2$ derivatives:
\be\label{eqn fio ests}
\rvts:L^2_r\left(X_R\right)\to L^2_{r-{\meff}^\sigma-\frac{q}2}\left(X_L\right).
\ee

It is natural to express the estimates for possibly degenerate FIO  in terms of possible losses 
relative to the optimal estimates.
Initially, our basic  assumptions is  that there is at least one $\sigma$ such that
(i) the double fibration condition \eqref{def DF general} is satisfied,  and 
(ii) there  is a known $\beta^\sigma\ge 0$ such that, for all $r\in\R$, 
\be\label{eqn basic k}
\rvts: L^2_r\left(X_R\right) \to L^2_{r-\meff^\sigma-\beta^\sigma}\left(X_L\right),
\ee
uniformly for $\vt\in T$, or at least for $\vt$ in some compact set containing any configurations that arise from the $E_i$ of interest.
\smallskip

\begin{remark}\label{rem uniform ests}
{A folk theorem in microlocal analysis is that the estimates for nondegenerate FIO  or even those covered by the 
corank $q$ scenario of \eqref{eqn fio ests}, which are all that we use in the 
concrete applications in this paper, are stable  under small perturbations of the amplitudes and phase functions in $C^N$ norm for $N$ sufficiently large. 
This is due to the finite number of integration by parts that are required in the various proofs for FIO
and the underlying oscillatory integral operators; 
see, e.g., \cite[Lem. 2.3]{GS94}. Thus, once one has a single value $\vt=\vt_0$ of the configuration parameter for which 
the generalized Radon transform $\rvts$ is nondegenerate or corank $q$, one is ensured that there is a neighborhood of 
$\vt_0$ for which this is true and for which \eqref{eqn fio ests} holds uniformly in $\vt$. 
See the comment at the end of Sec. \ref{sec tri pairs}.}
\end{remark}

Now  suppose that,  for $1\le i\le k$,   $E_i\subset X^i$ are compact sets.
Our goal is  to find  conditions on the $\hd(E_i)$ ensuring that 
$\Delta_\Phi\left(E_1,E_2,\dots,E_k\right)$
has nonempty interior in $T$.
For each $i$, fix an $s_i<\hd(E_i)$ 
and  a Frostman measure  $\mu_i$ on $E_i$ of finite $s_i$-energy;
translating energy into $L^2$-based Sobolev space norms,  $\mu_i\in L^2_{(s_i-d_i)/2}(X^i)$.
(See \cite{Mat95,Mat15} for further background.)
Define measures  
$$\mu_L := \mu_{i_1}\times \cdots \times \mu_{i_{|\sigma_L|}}\hbox{ on }X_L\hbox{ and }
\mu_R := \mu_{j_1}\times \cdots \times \mu_{j_{|\sigma_R|}}\hbox{ on }X_R,$$
and recall the following result from \cite{GIT20}:

\begin{proposition}\label{prop Sob}
For $1\le j\le k$, let $X^j$ be a $C^\infty$ manifold of dimension $d_j$, 
and suppose that  $u_j\in L^2_{r_j,\, comp}\left(X^j\right),\, 1\le j\le k$,
with each $r_j\le 0$. Then the tensor product
$u_1\otimes \cdots \otimes u_k$ belongs to 
$ L^2_{r,\, comp}\left(X^1\times \cdots \times X^k\right)$,
for $r=\sum_{j=1}^k r_j$.
\end{proposition}

From this it follows that
  $\mu_L\in L^2_{r_L}\left(X_L\right)$ and
 $\mu_R\in L^2_{r_R}\left(X_R\right)$, 
 where $r_L=\frac12\sum_{l=1}^{|\sigma_L|} \left(s_{i_l}-d_{i_l}\right)$ and
 $r_R=\frac12\sum_{l=1}^{|\sigma_R|} \left(s_{j_l}-d_{j_l}\right)$, resp.
 \medskip

As in \cite[Eqn. 2.6]{GIT20}, for any $\sigma\in\pk$, the configuration measure can be expressed as  
 \be\label{eqn pairing k}
 \nu(\vt)=\<\rvts\left(\mu_R\right),\mu_L\>,
 \ee
 which representation is justified {\it ex post facto} for $s_i$ in the admissible range.
 (See \cite[\S3.4]{GIT20} for the argument.)
Our basic assumption, that the boundedness \eqref{eqn basic k} holds for the $\sigma$ in question, then implies 
that $\rvts\left(\mu_R\right)\in L^2_{r_R-\meff^\sigma-\beta^\sigma}\left(X_L\right)$.
 Since $\mu_L\in L^2_{r_L}\left(X_L\right)$, the pairing in  \eqref{eqn pairing k} is bounded,
 and yields a continuous function of $\vt$ (by continuity of the integral), if
 \be\label{eqn first ineq}
 r_R-\meff^\sigma-\beta^\sigma + r_L \ge 0.
 \ee
 Noting that 
 $$r_L+r_R=\frac12\left[\left(\sum_{i=1}^k s_i\right) - d^{tot}\right],$$
 and using \eqref{eqn msigma}, we see that \eqref{eqn first ineq} holds iff
 \bast
 \sum_{i=1}^k s_i &\ge & d^{tot} +2\left(\meff^\sigma + \beta^\sigma\right)
 = d^{tot} + p - \min(d_L,d_R) + 2\beta^\sigma \\
 & = & \max(d_L,d_R) + p + 2\beta^\sigma.
  \east
 \medskip
 
  Optimizing over all nontrivial partitions $\sigma \in \pk$ leads to:
  \ms

\begin{theorem} \label{thm kpoint old}
{\bf Partition Optimization.} \cite[Thm. 5.2]{GIT20} \\  (i) With the notation and assumptions as above, define
\be\label{eqn po thresh}
s_\Phi=\min_\sigma \left[\,\max(d_L,d_R) + p + 2\beta^\sigma\, \right],
\ee
where the $\min$ is taken over those  $\sigma\in\pk$ for which both  the double fibration condition  
\eqref{def DF general} holds
and the uniform boundedness of the generalized Radon transforms $\rvt^\sigma$ with some loss of $\,\le \beta^\sigma$ derivatives 
\eqref{eqn basic k} hold.
Then, if $E_i\subset X^i$, $1\le i\le k$, are compact sets with $\sum_{i=1}^k\hd(E_i)>s_\Phi$, 
it follows that $\inter\left(\Delta_\Phi\left(E_1,E_2, \dots , E_k\right)\right)\ne\emptyset$.
\bs

(ii) In particular, if $X^1=\cdots=X^k=:X_0$, with $dim(X_0)=d$,
and $E\subset X_0$ is compact,
then $\inter\left(\Delta_\Phi\left(E\right)\right)\ne\emptyset$
if
\be\label{eqn equi k}
\hd(E)>\frac{1}{k}\left( \min_\sigma \max(d_L,d_R) + p \right),
\ee
where the minimum is taken over all
$\sigma\in\pk$ such that \eqref{def DF general} holds and the
canonical relation $C_{\vt}^\sigma$ is nondegenerate.  
\end{theorem}

The threshold for $\sum_{i=1}^k\hd(E_i)$ in \eqref{eqn po thresh} can be 
thought of as the 
{\it minimum} over all nontrivial partitions $\sigma$ of the thresholds 
determined by the \emph{maximum} microlocal loss (relative to the nondegenerate estimate)
over  all the points of $C_\vt^\sigma$. 
On general principle, one can (possibly) lower a {\it minimum of the maxima}  by replacing it with the \emph{maximum of the minima}, 
and in this setting it is  not hard to do this in practice.
The goal of this paper is to show that weakening the assumptions in the original partition optimization, by working 
{either locally on $\zvt$ or more generally}
microlocally on $N^*\zvt$,
can allow one to lower the needed threshold on $\hd(E)$,
or even to obtain a positive result when an application of the original version of partition optimization,
Thm. \ref{thm kpoint old}, 
would be  vacuous. 
\ms

In particular, in the context of Thm. \ref{thm kpoint old} (ii)  it is not necessary 
that \emph{any} of the canonical relations $C_\vt^\sigma$ be nondegenerate.
 Rather, working locally on $\zvt$, it is sufficient that, for every $x\in \zvt$ there is some neighborhood $U$ of $x$ in $\zvt$
and {\it some}  $\sigma\in\pk$ such that $C_\vt^\sigma$ is nondegenerate over $U$.
Even more generally, working microlocally, it suffices that
for every point $(x,\xi)\in N^*\zvt$,  
there exists some  $\sigma\in\pk$  and a conic neighborhood 
$\mathcal U$ of $(x,\xi)$ in $N^*\zvt$
such that  $C_\vt^\sigma$ is nondegenerate 
on $\mathcal U$ (or rather the image 
$\mathcal U^\sigma$ of $\mathcal U$ under the 
$\sigma$-separation of the variables to the left 
and right). Since a partition of unity subordinate to an open over of $\zvt$ is a special, $\xi$-independent case of a microlocal partition of unity
subordinate to a microlocal  cover of $N^*\zvt$,  the local version of the new approach is  a special case of the microlocal one.
However,  for clarity we state them separately:
\bs

\begin{theorem}\label{thm kpoint new local} {\bf Local Partition Optimization.} Suppose that there is a $\beta\ge 0$ such that, for
every point $x\in \zvt$ there exists a neighborhood $U$  and a partition $\sigma\in\pk$
for which  the generalized Radon transform $\rvt^\sigma$, localized to $U$,
satisfies both \eqref{def DF general} and  \eqref{eqn basic k} with a loss of at most $\beta$ derivatives, uniformly in $\vt$.

Then,
for $E\subset\R^d$ compact,   if
\be\label{eqn new equi k}
\hd(E)>\frac{1}{k}\left( 
\, \max(d_L,d_R) 
+p +2\beta\,  \right),
\ee
then $\inter\left(\Delta_\Phi\left(E\right)\right)\ne\emptyset$.
\end{theorem}
\bs

\begin{theorem}\label{thm kpoint new} {\bf Microlocal Partition Optimization.} Suppose  there exists a $\beta\ge 0$ such that,
for  every $(x,\xi)\in N^*\zvt$ there exist a conic neighborhood $\mathcal U$ and  partition $\sigma\in\pk$
for which  the generalized Radon transform $\rvt^\sigma$, microlocalized to $\mathcal U$,
satisfies both \eqref{def DF general} and  \eqref{eqn basic k} with a loss of at most $\beta$ derivatives, uniformly in $\vt$.
Then,
for $E\subset\R^d$ compact,  
$\inter\left(\Delta_\Phi\left(E\right)\right)\ne\emptyset$ if
\be\label{eqn new equi k}
\hd(E)>\frac{1}{k}\left( 
\, \max(d_L,d_R) 
+p +2\beta\,  \right).
\ee
\end{theorem}
\bs

Since spatial partitions of unity are special cases of microlocal ones,
the local theorem will  follow immediately from the microlocal one, 
which in turn is proven by  a straightforward refinement of the proof in \cite{GIT20}.
We start by forming a standard pseudodifferential partition of unity, 
$\sum Q_l(x,D)=I$, on $X$ subordinate to the open cover $\left\{\mathcal U_l\right\}$ of $N^*\zvt$,
supplemented by a $\mathcal U_0$ disjoint  from $N^*\zvt$ which completes the $\mathcal U$ 
to be a over of $T^*X\setminus 0$.
Each $Q_l\in \Psi^0_{cl}(X)$, and  together their principal symbols, $q_l(x,\xi)$, form a partition of unity on $T^*X\setminus 0$.
(For Thm. \ref{thm kpoint new local}, the $q_l$ are indpenendent of $\xi$.)
One can assume that this sum has at most $1+\left|\pk\right|$ terms.
We let $\sigma^l$ denote a partition such that $\rvt^\sigma$
satisfies \eqref{eqn basic k} with a loss of $\le\beta$ derivatives
on the conic support of $Q_l$.
The surface measure $\lambda_\vt$ from \eqref{def delta k} on $\zvt$ then decomposes as
$$\lambda_\vt = \sum_l Q_l\lambda_\vt,$$
leading to a similar decomposition of the generalized Radon transforms.
Hence, the identity \eqref{eqn pairing k} for the configuration measure can be replaced by
 \be\label{eqn new pairing k}
 \nu(\vt)=\sum_l \<\rvtsl\left(\mu_R\right),\mu_L\>.
 \ee
 By the analysis above, if $\hd(E)$ is greater than the threshold in \eqref{eqn new equi k},
 each of the terms in \eqref{eqn new pairing k} are continuous in $\vt$, finishing the proof.
 \ms
 
\begin{remark}\label{rem Hor}
We recall, for the proof  of Thms.  \ref{thm ratios}  below, that $\beta$ can be taken to be $r/2$
if the projection $\pi_L$  from each $\mathcal U$ drops rank by at most $r$  (see \cite{Hor71,Hor85}).
\end{remark}

\begin{remark}\label{rem nstartzt}
The conormal bundle of $\zvt$ is
 $$N^*\zvt=\left\{\left({x},D\Phi(\vec{x})^*(\tau)\right): {x}\in\zvt,\, \tau\in\R^p\setminus\mathbf 0\right\}.$$
 However, for the calculations needed to verify the microlocal condition in Thm. \ref{thm kpoint new}
 in each particular application,
  it is convenient to reorganize $N^*\zvt$ by grouping each pair $(x^i,\xi^i)\in T^* X^i$ together, and we define
 $$\widetilde{N^*\zvt}=\left\{\left(x^1,\xi^1; x^2, \xi^2; \dots; x^k,\xi^k\right): 
 (x^1,\dots,x^k;\xi^1,\dots,\xi^k)\in N^*\zvt\right\},$$
 and let $\pi_i$ denote the natural projection
 onto the $i$-th factor, $T^*X^i$.
 \end{remark}


\section{Areas of triangles}\label{sec triangles}

We now turn to results that require a microlocal approach, starting with the proofs of the various parts 
of Thm. \ref{thm triangles} concerning areas of triangles
generated by quadruples and quintuples of points in a planar set. 

\subsection{Pairs of areas of triangle in quadrilaterals}\label{sec tri pairs}

For part (i), let $\Phi:(\R^2)^4\to\R^2$ be
\bea\label{eqn phi i}
\Phi(x,y,z,w)&=&\big(\det[y-x,z-x],\det[z-x,w-x]\big)\nonumber\\
& &  \\
&=& \((y-x)\cdot(z-x)^\perp, (z-x)\cdot (w-x)^\perp\), \nonumber
\eea
where $\perp$ denotes rotation by $+\pi/2$, which is of course antisymmetric.
All of the entries in $D\Phi$ are $\perp$ of simpler expressions, and so in place of $D\Phi$ we work with
$$D\Phi^\perp:=
\left[\bmat
y-z & z-x & w-y & 0 \\
z-w & 0 & w-x & x-z
\emat\right],
$$
and we  denote the modified conormal bundle computed with $D\Phi^\perp$ by $\wtnzt^\perp$.

If $\hd(E)>5/3$, and $\mu$ is a Frostman measure for $s>5/3$, then since 
$\{(x,y,z)\in\R^6:  \det[y-x,z-x]=0\}$, the set of degenerate triangles, 
is an algebraic hypersurface, its Hausdorff dimension equals $5$.
 Hence,  $W_1:=\{(x,y,z,w): \det[y-x,z-x]=0\}$ has $\otimes^4 \mu$--measure 0 in $\R^8$, 
and without loss of generality we can assume that $4$-tuples we consider
lie in $\R^8\setminus W_1$; see \cite[Sec. 4.1]{GIT20} for related reasoning.
Thus, without loss of generality, we  can assume that
for each $\vt=(t_1,t_2)\in\R^2$, $\zvt=\Phi^{-1}(\vt)$ can be parametrized by $x,y,w\in\R^2$, 
with $z=x+\tz(x,y,w,\vt)\in\R^2$ then being the unique solution of
$$(x-y)^\perp\cdot (z-x)=t_1,\quad (w-x)^\perp\cdot(z-x)=t_2.$$
One can check that $|D\tz/Dy|\ne 0$ and $|D\tz/Dw|\ne 0$. 
\smallskip

Using the above one computes
\bea
\wtnzt^\perp&=& \Big\{\big(x,\tau_1(y-x-\tz)+\tau_2(x-w-\tz);\, y, \tau_1\tz;\nonumber\\
& &\quad  x+\tz,\tau_1(w-y)+\tau_2(w-x);\,  w,-\tau_2\tz\big):\, x,y,z\in\R^2,\, \tau\in\R^2\setminus\mathbf 0\Big\}.
\eea
From this we see that $D(x,\xi)/D(x,\tau)$ is always nonsingular. 

If $\tau_1\ne 0$, then $D(y,\eta)/D(y,w)$ is nonsingular, since $|D\tz/Dw|\ne 0$.
Thus, ordering the variables $x,y,z,w$, in order $1,2,3,4$,  partitioning them by $\sigma=(12|34)$ yields 
$C_\vt^\sigma$ which is  a local canonical graph on $\mathcal U_1=\{\tau_1\ne 0\}$.

On the other hand, if $\tau_2\ne 0$ then $D(w,\omega)/D(w,y)$ is nonsingular, since
$|D\tz/Dy|\ne \nolinebreak0$, so that using $\sigma=(14|23)$ gives $C_\vt^\sigma$ 
which is  a local canonical graph on $\mathcal U_2=\{\tau_2\ne 0\}$.

Together, $\mathcal U_1$ and $\mathcal U_2$ cover $\wtnzt^\perp$, and  $d_L=d_R=4$ for both partitions. 
{Picking any particular $\vt_0\in \Delta_\Phi(E)$ for which the above analysis applies, it also holds 
for $\vt$ close to $\vt_0$, and the estimates for resulting nondegenerate generalized Radon 
transforms $\rvts$ are locally uniform in $\vt$; cf. Remark \ref{rem uniform ests}.}
Thus, Thm. \ref{thm kpoint new} applies with $\beta=0$. 
It follows that if $\hd(E)>\frac14(4+2+0)=3/2$, $\Delta_\Phi(E)$ has nonempty interior in $\R^2$.


\subsection{Triples of areas of triangles in quadrilaterals}\label{sec tri triples}

To prove  Thm. \ref{thm triangles}(ii)
we modify the considerations of the previous section as follows.
Let $\Phi:(\R^2)^4\to\R^3$ be
\bea\label{eqn phi ii}
\Phi(x,y,z,w)&=&\big(\det[y-x,z-x],\det[z-x,w-x],\det[y-x,w-x]\big)\nonumber\\
& &  \\
&=& \((y-x)\cdot(z-x)^\perp, (z-x)\cdot (w-x)^\perp, (w-x)\cdot (y-x)^\perp\). \nonumber
\eea
As before,  in place of $D\Phi$ we work with
$$D\Phi^\perp:=
\left[\bmat
y-z & z-x & x-y & 0 \\
z-w & 0 & w-x & x-z \\
w-y & x-w & 0 & y-x & 
\emat\right],
$$
and denote the modified conormal bundle computed with $D\Phi^\perp$ by $\wtnzt^\perp$.

For $\vt=(t_1,t_2,t_3)\in\R^3$, $\zvt=\Phi^{-1}(\vt)$ is determined by
$$(y-x)^\perp\cdot  (z-x)=-t_1,\, (w-x)^\perp\cdot (z-x)=t_2,\, (y-x)^\perp\cdot (w-x)= t_3.$$
Solving the last equation first, we can solve for $w$ with one degree of freedom:
\bast
w&=\,\,& x+t_3\frac{(y-x)^\perp}{|y-x|}+s(y-x),\, s\in\R\\
&=:& x+ \tw(x,y,s;t_3),
\east
so that $w-x=\tw$. Then, as in the previous section, 
without loss of generality we can assume that $\det[y-x,z-x]\ne0$ and so
one can solve uniquely for $z$, incorporating the 
dependence of $w$ on $s$:
$$z=x+\tz(x,y,s;\vt) \implies z-x=\tz.$$
Note that $\partial_s\tw=y-x$ and, as in the previous section, $|D\tz/Dy|\ne 0$, 
\smallskip

We can  parametrize the conormal bundle as

\bea
\nonumber
\wtnzt^\perp&=&\Big\{\big( x,\tau_1(y-x-\tz)+\tau_2(\tz-\tw)+\tau_3(\tw-y); \\
& & \qquad y,\tau_1\tz-\tau_3\tw;\, x+\tz,\tau_1(x-y)+\tau_2\tw;\\
& & \qquad\quad x+\tw,-\tau_2\tz+\tau_3(y-x)\big): \, x,y\in\R^2, s\in\R, \tau\in\R^3\setminus\mathbf 0\Big\}.
\nonumber
\eea

Note that the differential of $(x,\xi)$ with respect to $x$ and any two of the three $\tau_j$ is nonsingular.
(Here we can assume that $y-x, z-w$ and $w-y$ are in general position, i.e., any two are linearly 
independent, which excludes a variety  $W_2\subset\R^8$ of dimension 5.) 
This leaves $y, s$ and the remaining $\tau_j$ variable to use for another one of the three remaining projections.

Since $\partial_s\tw=y-x$, one sees that  $D(y,\eta)/D(y,s,\tau_1)$ is nonsingular if $\tau_3\ne0$,
while $D(y,\eta)/D(y,s,\tau_3)$ is nonsingular if $\tau_1\ne0$.
Hence, $C_\vt^{(12|34)}$ is a local canonical graph on 
$\mathcal U_1=\{\tau_1\ne0\hbox{ or } \tau_3\ne0\}$.
\smallskip

Combining $\partial_s\tw=y-x$ with $|D\tz/Dy|\ne 0$, one sees that $D(z,\zeta)/D(y,s,\tau_1)$ is 
nonsingular if $\tau_2\ne 0$, so that $C_\vt^{(13|24)}$ is a local canonical graph on 
$\mathcal U_2=\{\tau_2\ne 0\}$.

Since $\mathcal U_1,\, \mathcal U_2$ form an open cover of $\wtnzt^\perp$,
and $d_L=d_R=4,\, \beta=0$ for all of those partitions,
we can apply  Thm. \ref{thm kpoint new}, obtaining that 
 if $\hd(E)>\frac14(4+3+0)=7/4$, $\Delta_\Phi(E)$ has nonempty interior in $\R^3$.


\subsection{Triples of areas of triangles in pentagons}\label{sec tri pentagon}

For the proof of  Thm. \ref{thm triangles} (iii)
we modify the setup for parts (i) and (ii)  as follows.
Define $\Phi:(\R^2)^5\to\R^3$, recording the areas of the three adjacent triangles pinned at $x$, 
by
\bea\label{eqn phi ii}
\Phi(x,y,z,w,u)&=&\big(\det[y-x,z-x],\det[z-x,w-x],\det[w-x,u-x]\big)\nonumber\\
& &  \\
&=& \((y-x)\cdot(z-x)^\perp, (z-x)\cdot (w-x)^\perp, (w-x)\cdot (u-x)^\perp\). \nonumber
\eea
As before,  in place of $D\Phi$ we work with
$$D\Phi^\perp:=
\left[\bmat
y-z & z-x & x-y & 0 & 0 \\
z-w & 0 & w-x & x-z & 0\\
w-u & 0 & 0 & u-x & x-w
\emat\right],
$$
and denote the modified conormal bundle computed with $D\Phi^\perp$ by $\wtnzt^\perp$:

\bea\label{eqn pentagon general}\nonumber
\wtnzt^\perp &=& \Big\{\big(x,\tau_1(y-z)+\tau_2(z-w)+\tau_3(w-u);\, y,\tau_1(z-x);\\
& & \quad\,z,\tau_1(x-y)+\tau_2(w-x);\, w,\tau_2(x-z)+\tau_3(u-x);\\
& & \quad\,u,\tau_3(x-w)\big):\, (x,y,z,w,u)\in\zvt,\, \tau\in\R^3\setminus\mathbf 0\Big\}.
\nonumber
\eea

The linear coordinates $\tau_1,\tau_2,\tau_3$ on the fibers are intrinsically defined 
(given that $\Phi$ has been fixed), independent of what coordinates we pick  on the $7$-dimensional
base $\zvt$. We claim that on the open conic sets
$\mathcal U_j=\{ \tau_j\ne0\}\subset \wtnzt^\perp,\, j=1,2,3$, which  form  a microlocal cover, 
the partitions $\sigma=(14|235),\, (13|245),\, (13|245)$, resp., give canonical relations $C_\vt^\sigma$ 
which are nondegenerate. (Note that the partitions used on $\mathcal U_2$ and $\mathcal U_3$ are the 
same, but we have to treat $\mathcal U_2$ and $\mathcal U_3$ separately.)
Thus, Thm. \ref{thm kpoint new} implies that $\Delta_\Phi(E)$  
has nonempty interior for for $E\subset \R^2$ with $\hd(E)>\frac15\(\max(4,6)+3+0\)=\frac95$,
proving Thm. \ref{thm triangles}(iii).
\smallskip

To prove the claim above, we use two different coordinate parametrizations of $\zvt$; the first 
is useful for establishing the claim on $\mathcal U_1$ and $\mathcal U_2$, 
and the second for $\mathcal U_3$.
For the first, we parametrize $\zvt$  by $x,y,w\in\R^2$ and $ s\in\R$ by 
\smallskip

(i) solving the $2\times 2$ system for $z$,
$$(y-x)^\perp\cdot  (z-x)=-t_1,\, (w-x)^\perp\cdot (z-x)=t_2,$$
obtaining, for $(x,y,w)$ in  general position (in the complement of a hypersurface), a unique solution
$z=x+\tz(x,y,w,\vt)$; and 
\smallskip

(ii) for $w\ne x$, solving 
$ (w-x)^\perp\cdot (u-x) = -t_3$ for $u$ with one degree of freedom, 
$$u-x=-t_3\frac{(w-x)^\perp}{|w-x|} + s\cdot(w-x)=:\tilde{u}(x,w,s;t_3),\, s\in\R.$$
Note that 
\be\label{eqn pentagon Dyw}
\left| \frac{D(y-x-\tz)}{Dy} \right| \ne 0,\quad \left| \frac{D(w-x-\tz)}{Dw} \right| \ne 0;
\ee
the first follows since the differential maps $(y-x)\cdot\partial_y \to (y-x)\cdot\partial_y$ and 
$(y-x)^\perp\cdot \partial_y \to c_{y,w,\vt}(y-x)^\perp\cdot \partial_y +\dots$,
and the second is similar. We also have $\partial_s\tu=w-x\ne 0$.
\smallskip

Adapting \eqref{eqn pentagon general} to this parametrization of $\zvt$, the conormal bundle of $\zvt$ is
parametrized

\bea\label{eqn pentagon A}
\nonumber
\wtnzt^\perp&=&\Big\{\big(x,\tau_1(y-x-\tz)+\tau_2(x-w+\tz)+\tau_3(w-x+\tu);\, y, \tau_1\tz;  \\
& & \quad\,  x+\tz,\tau_1(x-y)+\tau_2(w-x);\, w,-\tau_2\tz+\tau_3\tu;\, x+\tu, \tau_3(x-w)\big) \\
& & \qquad  :\, x,y,w\in\R^2,s\in\R,\tau\in\R^3
\setminus\mathbf 0
\Big\}.
\nonumber
\eea

On $\{\tau_1\ne0\}$, we can use the $\tau_1$ term in the expression for $\xi$
in \eqref{eqn pentagon A} together with \eqref{eqn pentagon Dyw} to obtain $|D(x,\xi)/D(x,y)|\ne 0$,
while $D(w,\omega)/D(w,\tau_2,\tau_3)|\ne0$ since $\tz,\, \tu$ are generically linearly independent.
Hence, $\pi_1\times\pi_4:C_\vt^{(14|235)}\to T^*\R^4$ is a submersion.

On $\{\tau_2\ne0\}$, $|D(x,\xi)/D(x,y)|\ne 0$ using   \eqref{eqn pentagon Dyw}
with the $\tau_2$ term in the expression for $\xi$ in \eqref{eqn pentagon A},
while $|D(z,\zeta)/D(y,\tau_1,\tau_2)|\ne0$ from \eqref{eqn pentagon Dyw} 
and the generic linear independence of $x-y,\, w-x$. 
Hence, $\pi_1\times\pi_3:C_\vt^{(13|245)}\to T^*\R^4$ is a submersion.
\medskip

To deal with $\mathcal U_3=\{\tau_3\ne0\}$, we change the parametrization to  
$x,z,u\in\R^2$ and $s'\in\R$ by 

(i) solving the $2\times 2$ system for $w$,
$$(z-x)^\perp\cdot  (w-x)=-t_2,\, (u-x)^\perp\cdot (w-x)=t_3,$$
obtaining, for $(x,z,w)$ in  general position a unique solution
$w=x+\tw(x,z,u;\vt)$, with $\left|D(u-x-\tw)/Du\right|\ne 0$; and 
\smallskip

(ii) for  $z\ne x$, solving 
$ (z-x)^\perp\cdot (y-x) = t_1$ for $y$ with one degree of freedom, 
$$y-x=t_1\frac{(z-x)^\perp}{|z-x|} + s'\cdot(z-x)=: \tilde{y}(x,z,s';\vt),\, s'\in\R.$$

With respect to these coordinates, the analogue of \eqref{eqn pentagon A} is

\bea\label{eqn pentagon B}
\nonumber
\wtnzt^\perp&=&\Big\{\big(x,\tau_1(x-z+\ty)+\tau_2(z-x-\tw)+\tau_3(x-u+\tw);\, x+\ty, \tau_1(z-x);  \\
& & \quad\,  z,-\tau_1\ty+\tau_2\tw;\, x+\tw,\tau_2(x-z)+\tau_3(u-x);\, u, -\tau_3\tw \big) \\
& & \qquad  :\, x,z,u\in\R^2,s'\in\R,\tau\in\R^3
\setminus\mathbf 0
\Big\}.
\nonumber
\eea
Arguing as above, one sees that, using the $\tau_3$ term in $\xi$, for $\tau_3\ne 0$ we have
$\left|D(x,\xi)/D(x,u)\right|\ne0$, and $\left|D(z,\zeta)/D(z,\tau_1,\tau_2)\right|\ne0$, which shows
that $\pi_1\times\pi_3$ is a submersion. 
Hence, $C_\vt^{(13|245)}$ is nondegenerate on $\mathcal U_3$, and this finishes the proof of 
Thm. \ref{thm triangles}(iii).


\section{Pairs of ratios of  distances}\label{sec ratios}

We now discuss and prove Thm. \ref{thm ratios}: 
if $E\subset\R^d$ with $\hd(E)>(3d+1)/4$, then 
\be\label{thm 1.2 eqn}
\hbox{Int }\left\{\left(\frac{|x-y|}{|z-w|},\, \frac{|x-w|}{|z-y|}\right)\in\R^2:
x,y,z,w\in E,\, z\ne w,\, z\ne y \right\}\ne\emptyset.
\ee
See the left 4-tuple in Figure \ref{fig one} below.

To put this in context, we discuss sets of ratios of distances more generally,
recalling  related previous results and then  proving some variations.
It follows immediately from the result  of Mattila and Sj\"olin \cite{MS99} that
if $E\subset \R^d,\, \hd(E)>(d+1)/2$, then
\be\label{eqn trivial}
\hbox{Int }\left\{\frac{|x-y|}{|z-w|}\, :\, x,y,z,w\in E,\, z\ne w\right\}\ne\emptyset,
\ee
since both the numerators and denominators range over sets containing an interval.
\ms

Later, under the higher dimensional threshold $\hd(E)>(d+2)/2$, 
Peres and Schlag \cite[p. 248]{PSch00}  proved a stronger, pinned version of \cite{MS99}:
there exists an $x\in E$ such that the pinned distance set,
$$\Delta^x(E):=\left\{\,|x-z|\, :\, z\in E\, \right\}$$
contains an interval, i.e., has nonempty interior.  
Using such a pin point $x$, and then fixing any $y\in E,\, y\ne x$,
it  follows immediately that for $\hd(E)>(d+2)/2$,

\be\label{result ps}
\left(\exists\, x\in E\right)\, \left(\forall\, y\in E,\, y\ne x\right)\, \hbox{Int }\left\{\,\frac{|x-z|}{|x-y|}\, :\, z\in E\, \right\} \ne\emptyset.
\ee

By weakening the statement in \eqref{result ps},
we can lower the threshold:

\begin{theorem}\label{thm yz}
If $E\subset\R^d,\, d\ge 2$, is compact and $\hd(E)>(d+1)/2$, then
\be\label{result ps modified}
\left(\exists\, x\in E\right)\,\hbox{ \emph{s.t.} }  \inter\left\{\,\frac{|x-z|}{|x-y|}\, :\, y,\, z\in E,\, y\ne x\, \right\} \ne\emptyset.
\ee
\end{theorem}

\begin{remark}\label{rem git thm1.4}
This improves upon our result \cite[Thm. 1.4]{GIT20}, which established \eqref{result ps modified}
for $\hd(E)>(2d+1)/3$ using the original partition optimization.
\end{remark}

\begin{remark}\label{rem bio23} The elementary proof we give, combining known results, 
follows closely an argument 
in a recent preprint of Borges, Iosevich and Ou \cite[Sec. 2]{BIO23}.
\end{remark}

To prove Theorem \ref{thm yz}, first recall  another result of Peres and Schlag  \cite[Cor. 8.4]{PSch00}
(see also \cite{ITU16,IL19}):  if $\hd(E)>(d+1)/2$, there exists $x\in E$ such that the pinned distance set
has positive Lebesgue measure, $\left|\Delta^x(E)\right|_1>0$.
Apply this twice
 for one such  pin point, $x$,  with  variables $z$ and $y$, resp., in the definition of $\Delta^x(E)$.
For each of these, apply $\log$ (a bi-Lipschitz map on compact subsets of $\R_+$)  to obtain
$$\left| \{ \, \log |x-z|\, :\, z\in E,\, z\ne x\,\}\right|_1 >0,\quad \left| \{ \, \log |x-y|\, :\, y\in E,\, y\ne x\,\}\right|_1 >0.$$
Now invoking the classical theorem of Steinhaus \cite{Steinhaus} on difference sets, one has that
$$\inter \left\{\,  \log |x-z| - \log |x-y|:  \,  y,z\in E,\,\, y\ne x,\, z\ne x \, \right\}\ne\emptyset.$$
Applying the exponential map (a diffeomorphism $\R \to \R_+$) preserves nonempty interior, 
yielding   \eqref{result ps modified} and finishing the proof of Theorem \ref{thm yz}.
\ms

These results are motivation to study other configuration sets  constructed from ratios of distances. 
For example, one such is the following: 
If $E\subset \R^d,\, d\ge 2$, is  compact with $\hd(E)>(d+1)/2$, then
\be\label{eqn zywy}
\left(\exists\, x\in E\right)\, \left(\forall\, y\in E,\, y\ne x\right)\,\hbox{Int }\left\{\left(\frac{|x-z|}{|x-y|},\, 
\frac{|x-w|}{|x-y|}\right)\in\R^2: 
z,w\in E\,   \right\}\ne\emptyset.
\ee
(See the right 4-tuple in Figure \ref{fig one} below.)
However, this follows immediately from \eqref{result ps} since, with $x$ and $y$ fixed,  the set in \eqref{eqn zywy}
\atxt{equals}  the Cartesian product of the set in \eqref{result ps} with itself \atxt{and hence has nonempty interior.}
\ms

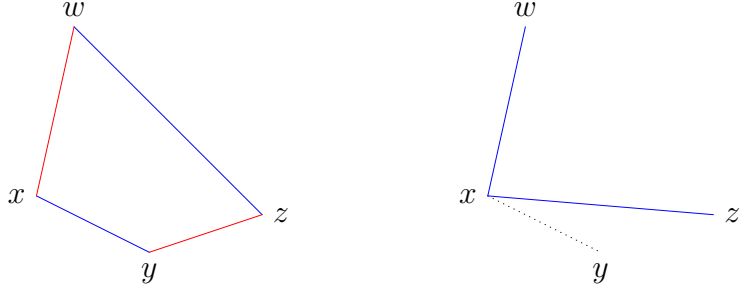
\begin{figure}\label{fig one}
    \centering
    \begin{tikzpicture}
    
        \tkzDefPoint(0,0.75){x}
        \tkzDefPoint(1.5,0){y}
        \tkzDefPoint(3,0.5){z}
         \tkzDefPoint(0.5,3){w}
        \tkzDrawSegments[color = blue](x,y)
           \tkzDrawSegments[color = blue](z,w)
             \tkzDrawSegments[color = red](z,y)
           \tkzDrawSegments[color = red](x,w)
        
        \tkzLabelPoint[left](x){$x$}
        \tkzLabelPoint[below](y){$y$}
        \tkzLabelPoint[right](z){$z$}
        \tkzLabelPoint[above](w){$w$}
        
        \tkzDefPoint(6,0.75){x'}
        \tkzDefPoint(7.5,0){y'}
        \tkzDefPoint(9,0.5){z'}
         \tkzDefPoint(6.5,3){w'}
        \tkzDrawSegments[dotted](x',y')
           \tkzDrawSegments[color = blue](z',x')
             \tkzDrawSegments[color = blue](w',x')
        
        \tkzLabelPoint[left](x'){$x$}
        \tkzLabelPoint[below](y'){$y$}
        \tkzLabelPoint[right](z'){$z$}
        \tkzLabelPoint[above](w'){$w$}
        
    \end{tikzpicture}
    \caption{Two pairs of ratios of distances of points in a 4-tuple. 
    Figure on left corresponds to \eqref{thm 1.2 eqn} from Thm. \ref{thm ratios},
    with ratio of the two blue side lengths forming the first coordinate of $\Phi(x,y,z,w)$ 
    and the ratio of the two red side lengths being the second coordinate.
    Figure on right illustrates \eqref{eqn zywy} for the same 4 points, 
    with the two coordinates of $\Phi(x,y,z,w)$ being the ratios of the lengths of the blue segments to the dotted one.}
    \label{fig one}
\end{figure}

However, the set  that is the subject of  Thm. \ref{thm ratios}, i.e.,  the one in \eqref{thm 1.2 eqn}, 
does not seem to be amenable to such an approach, 
and we prove it using microlocal partition optimization, Thm. \ref{thm kpoint new},with $k=4,\, p=2$,
which results in requiring the higher dimensional threshold, $\hd(E)>(3d+1)/4$.

To start the proof,  the configuration function $\Phi:\left(\R^d\right)^4\to \R_+^2$, 
\be\label{eqn ratio phi}
\Phi(x,y,z,w):=\left(\frac{|x-y|}{|z-w|},\, \frac{|x-w|}{|z-y|}\right),\,
\ee
defines, for each $\vt=(t_1,t_2)\in\left(\R_+\right)^2$, an  incidence relation $\zvt=\Phi^{-1}(\vt)$.
Defining $F=(F_1,F_2)$, with
$$
F_1=|x-y|-t_1|z-w|,\quad F_2=|x-w|-t_2|z-y|,$$
so that $\zvt=F^{-1}(\mathbf 0)$, 
a calculation shows that $\zvt$  can be parametrized by
 \bea\label{eqn ratio zt}\nonumber
\zvt&=&\Big\{(x,y,z,w)=(x,x-t_1r\omega^1,z,z-r\omega^2)\in \R^{4d}\\
& & \qquad :\, x,z\in\R^d,\, x\ne z,\, r>0,\, \left(\omega^1,\omega^2\right)\in S_{x,z,r} \, \Big\},
\eea
where
\bast
S_{x,z,r}  &=& \big\{\, (\omega^1,\omega^2)\in \sd\times\sd \\
& & \qquad :\, 2(z-x)\cdot(t_1t_2^22\omega^1-\omega^2)=(1-t_2)^2|x-z|^2+(1-t_1^2t_2^2)r^2\, \big\}.
\east
One further calculates
\bast
DF&=&\left[
\begin{matrix}
\frac{x-y}{|x-y|} & -\frac{x-y}{|x-y|} & -t_1\frac{z-w}{|z-w|} & t_1\frac{z-w}{|z-w|}  \\
\frac{x-w}{|x-w|} & t_2\frac{z-y}{|z-y|} & -t_2\frac{z-y}{|z-y|} & -\frac{x-w}{|x-w|}
\end{matrix}
\right] \\
& & \\
& & \\
&=& \left[
\begin{matrix}
\omega^1 & -\omega^1 & -t_1\omega^2 & t_1\omega^2 \\
\frac{x-z+r\omega^2}{|x-z+r\omega^2|} & t_2\frac{z-x+t_1r\omega^1}{|z-x+t_1r\omega^1|} &
-t_2\frac{z-x+t_1r\omega^1}{|z-x+t_1r\omega^1|} & -\frac{x-z+r\omega^2}{|x-z+r\omega^2|} 
\end{matrix}
\right],
\east
\smallskip

\noindent where the second representation is $DF$  evaluated at points of $\zvt$ in 
terms of the parametrization \eqref{eqn ratio zt}. As in Remark \ref{rem nstartzt}, we let
$$\wtnzt=\left\{\left(x,\tau^tD_xF;\, y,\tau^tD_yF;\, z,\tau^tD_zF;\, w,\tau^tD_wF\right)\, 
:\, (x,y,z,w)\in\zvt,\, \tau\in\R^2\setminus 0\, \right\}$$
be the conormal bundle of the incidence relation, with the $(x,\xi),\, (y,\eta),\, (z,\zeta),\, (w,\nu)$ variables
separated.

Denote the  projections from $\widetilde{N^*\zvt}$ onto its four $T^*\R^d$ factors by $\pi_x,\pi_y,\pi_z,\pi_w$,
respectively,
each of which is a function of $x,z,r,\omega,\tau$ with values in $T^*\R^d$. 
Corresponding to the choice of partition $\sigma=(13\, |\, 24)$, by the above  one has
$$(\pi_x\times \pi_z)(x,z,r,\omega,\tau)=
\Big(\, x,\tau_1\omega^1+\tau_2\frac{x-z+r\omega^2}{|x-z+r\omega^2|};\,
z, -\tau_2t_1\omega^2-\tau_2 t_2\frac{z-x+t_1r\omega^1}{|z-x+t_1r\omega^1|}\, \Big).$$
From this one can calculate that, 
for $\tau_1\ne 0$ and $x-z+t_1r\omega^1\ne 0$, 
$$\hbox{Rank } \frac{D(x,\xi,z,\zeta)}{D(x,z,\omega^1,\tau)}=3d+1,$$
and similarly, for $\tau_2\ne 0$ and $x-z+r\omega^2\ne 0$, 
$$\hbox{Rank } \frac{D(x,\xi,z,\zeta)}{D(x,z,\omega^2,\tau)}=3d+1.$$
Note that $W:=\{x-z+t_1r\omega^1= 0\}\cup\{x-z+r\omega^2= 0\}$ is a variety of dimension $3d$.
Since our dimensional threshold will be $\hd(E)>(3d+1)/4$, and thus $\hd(E\times E\times E\times E)>3d+1 > 3d$,
this exceptional set is irrelevant for the analysis. See, e.g., \cite[Sec. 4]{GIT20} for several instances of this type of 
reasoning. Thus, if over $\zvt\cap\left(\R^{4d}\setminus W\right)$ 
we let $\mathcal U_j=\{\tau_j\ne0\}\subset \wtnzt,\, j=1,2$,
then $\{ \mathcal U_1,\, \mathcal U_2\}$ forms an open cover of $C^\sigma_\vt$ on which 
$D(\pi_x\times \pi_z)$ drops rank by $\le d-1$.
Hence, by the FIO estimates (cf. Remark \ref{rem Hor}), the associated  $\rvts$ loses at most $\beta^\sigma=(d-1)/2$ derivatives.
Applying Theorem \ref{thm kpoint new}  with codimension $p=2$ and loss $\beta=(d-1)/2$ then yields
that $\Delta_\Phi(E)$, the set in \eqref{thm 1.2 eqn}, has nonempty interior if
$\hd(E)>\frac14(2d+2+(d-1))=(3d+1)/4$, finishing the proof of Thm. \ref{thm ratios}.


\section{Congruence classes of triangles in $\R^d$}\label{sec PRA}

We now   show how  the result of Palsson and Romero-Acosta \cite{PRA21} follows easily from the 
microlocal approach taken here. Let $\Phi:\( \R^d \)^3\to\R^3$,
$$\Phi(x,y,z)=\(\,|x-y|,\, |x-z|,\, |y-z|\, \),$$
so that $\Delta_\Phi(E)$ is the set of vectors of side lengths of triangles generated by
the points in $E$ and thus, 
modulo permutations, the set of congruence classes of triangles in $E$.
Letting $\widehat{x}=x/|x|$, one computes
$$D\Phi=
\[\bmat
\widehat{x-y} & -\widehat{x-y} & 0 \\
\widehat{x-z} & 0 & -\widehat{x-z} \\
0 & \widehat{y-z} & -\widehat{y-z}
\emat\].$$
Furthermore, for $\vt=(t_1,t_2,t_3)\in\R_+^3$,
we can parametrize $\zvt$ as follows: We first take $x\in\R^d$ to be arbitrary, 
and then  $y=x-t_1\omega$,  with $\omega\in\sd$ arbitrary.
If one writes $z=x-t_2\tom$ for some $\tom\in\sd$, 
then one computes that $|y-z|=t_3$ iff $\tom\cdot\om=\frac{t_3-(t_1^2+t_2^2)}{2t_1t_2}$.
For $\vt$ in the complement of a lower dimensional variety,
$$S_{\vt,\om}:=\left\{ \tom\in\sd: \tom\cdot\om=\frac{t_3-(t_1^2+t_2^2)}{2t_1t_2}\right\}$$
is a smooth $(d-2)$-surface in $\sd$ (possibly empty), and 
$$\zvt=\left\{ (x,x-t_1\om,x-t_2\tom)\, :\, x\in\R^d, \om\in\sd, \tom\in S_{\vt,\om}
\right\}.$$
Applying $D\Phi^*$ to $\tau\in\R^3$ at these points, we obtain
\bea\nonumber
\wtnzt&=&\Big\{\big(x,\, \tau_1\om+\tau_2\tom;  
x-t_1\om,\, -\(\tau_1-(t_1/t_3)\tau_3\)\om +(t_2/t_3)\tau_3\tom; \\
& & \qquad \qquad\qquad\qquad\,\, \, 
x-t_2\tom,\, (t_1/t_3)\tau_3\om -\(\tau_2+(t_2/t_3)\tau_3\)\tom\big) \\
& & \qquad \qquad\qquad\qquad\,\, \,\,  : \, x\in\R^d,\, \om\in\sd,\, \tom\in S_{\vt,\om},\, 
\tau\in\R^3\setminus\mathbf 0\Big\}.
\nonumber
\eea

Let $i_\om:T_\om\sd\hookrightarrow T_\om\R^d$ and 
$\ti:T_{\tom} S_{\vt,\om}\hookrightarrow T_{\tom}\R^d$
{ be the inclusions of the tangent spaces}, and note that for generic $\vt$, 
\be\label{eqn tri span}
span\left\{T_\om\sd,\, \om\right\}=T_\om\R^d\hbox{ and } 
span\left\{T_{\tom} S_{\vt,\om},\,\om,\,\tom\right\}
=T_{\tom}\R^d.
\ee
Denoting the projections into the $(x,\xi),\, (y,\eta)$ and $(z,\zeta)$ variables by $\pi_j,\, j=1,2,3$, resp., we 
calculate their Jacobians with respect to $(x,\om,\tom,\tau_1,\tau_2,\tau_3)$; 
to avoid clutter, we indicate unneeded terms by $*$:

$$D\pi_1=
\[\bmat
I_d & 0 & 0 & 0 & 0 & 0 \\
0 & \tau_1i_\om & \tau_2\ti & \om & \tom & 0
\emat\],$$

$$D\pi_2=
\[\bmat
I_d & * & * & * & * & * \\
0 & -\(\tau_1-(t_1/t_3)\tau_3\)i_\om & (t_2/t_3)\tau_3\ti & -\om & 0 & (1/t_3)(t_2\tom+t_1\om)
\emat\],$$
and
$$D\pi_3=
\[\bmat
I_d & * & * & * & * & * \\
0 & (t_1/t_3)\tau_3i_\om &  -\(\tau_2+(t_2/t_3)\tau_3\) \ti & 0 & -\tom & (1/t_3)(t_1\om -t_2/\tom)
\emat\].$$
\smallskip 

Examining the column spaces of these and using \eqref{eqn tri span}, one sees that $D\pi_1$ is surjective except on the line 
$$L_1:=\{\tau_1=\tau_2=0\};$$
$D\pi_2$ is surjective except on 
$$L_2:=\{\tau_1+(t_1/t_3)\tau_3=\tau_3=0\}=\{\tau_1=\tau_3=0\};$$ 
and $D\pi_3$ is surjective except  on 
$$L_3:=\{\tau_3=\tau_2+(t_2/t_3)\tau_3=0\}=\{\tau_2=\tau_3=0\}.$$
Furthermore, for each $j$,  the image $\pi_j(L_j)$ lies in the $\mathbf0$-section of $T^*\R^d$,
which causes problems with the standard theory of FIO.
Since this holds for every partition, 
the original partition optimization of Thm.  \ref{thm kpoint old} is inapplicable. 
Furthermore, since the lines of degeneracy exist above all points of $\zvt$,
the merely local version, Thm.  \ref{thm kpoint new local}, also does not suffice, 
so that one needs  the full strength of the microlocal version.
\smallskip

Setting $\mathcal U_j=\{\tau\in\R^3\setminus L_j\}\subset \wtnzt$, 
it follows that $\{\mathcal U_1,\mathcal U_2,\mathcal U_3\}$ is a conic open cover of $\wtnzt$ 
on which the partitions $\sigma=(1|23),\, (2|13),\, (3|12)$, resp., result in canonical relations
$C_\vt^\sigma\subset \left(T^*\R^d\setminus \mathbf 0\right)\times \left(T^*\R^{2d}\setminus \mathbf0\right)$ 
which are nondegenerate, so that Thm. \ref{thm kpoint new} applies with 
$k=3,\, p=3, \max(d_L,d_R)=2d$ and $\beta=0$. Hence, for $E\subset\R^d$ with
$$\hd(E)>\frac13\(2d+3+0\)=(2d+3)/3,$$
$\Delta_\Phi(E)$ has nonempty interior in $\R^3$, reproving the main result of \cite{PRA21}.
\medskip

As a final comment, we remark that in their  recent  preprint \cite{PRA22}, Palsson and Romero-Acosta 
have extended the results of \cite{PRA21} to $(k-1)$-simplices in $\R^d$, $k\ge 4$, 
for some thresholds depending on $d$ and $k$.
Calculations along the lines of those used above indicate that some of the conditions 
 required to apply Thm. \ref{thm kpoint new} for these higher values of $k$ appear to 
 fail. It would be interesting to see whether  further
microlocal analysis of the problem can be used to obtain results for these higher dimensional simplices.

\section{Similarity classes of triangles in $\R^d$}\label{sec similarity}

We conclude with the proof of Thm. \ref{thm similarity}. 
This is a  variation on the  result of the previous section,
with congruence of triangles replaced by similarity and the resulting threshold  lowered by $1/3$;
however,  it can also be viewed as concerning ratios of distances, in the spirit of  Sec. \ref{sec ratios}.
Just as \eqref{eqn pra} encodes   the congruence classes of triangles with vertices in $E$
(up to  permutations of $x,y,z$, which do not affect the nonempty interior statement),
\eqref{eqn sim} encodes the similarity classes of triangles with vertices in $E$.
To make this more explicit, recall that the projective plane $\rp$ is covered by three coordinate charts, 
$$\mathcal V_\alpha:=\big\{ \, [\,A\,:\, B\,:\, C\,]
\hbox{ s.t. }A,B,C\in\R,\, \alpha\ne 0\, \big\},\quad \alpha=A,\, B\hbox{ or } C.$$
In particular, $\mathcal V_A=\big\{\, [1\, :\, u \, :\, v\, ]\hbox{ s.t. } (u,v)\in\R^2\, \big\}$, and to prove 
Thm. \ref{thm similarity} it suffices to show that  the configuration function
$$\Psi(x,y,z):=\left(\frac{|x-z|}{|x-y|}\, ,\, \frac{|y-z|}{|x-y|}\right) $$
satisfies $\hbox{Int }(\Delta_\Psi(E))\ne\emptyset$ for $\hd(E)>(2d+2)/3$. 
\ms

For $\vt=(t_1,t_2)\in \R_+^2$, we parametrize $\zvt$ as follows: For $x\in \R^d$,  write $y=x+r\omega$ 
with $r>0,\,  \omega\in\sd$, so that $|x-y|=r$.  In order for $(x,y,z)\in\zvt$, $z$ must be of the form $z=x+t_1r\tom$ 
for some $\tom\in\sd$, ensuring that $|x-z|=t_1r$; the further constraint that $|y-z|=t_2|x-y|$ then becomes 
$$\omega\cdot\tom=(2t_1)^{-1}\left( 1+t_1^2-t_2^2\right):=a(\vt).$$
Thus $\tom=a(\vt)\om+b(\vt)\ttom$ for some $\ttom\in\sd\cap\,  \om^\perp$, 
where $b(\vt)=\left(1-a^2(\vt)\right)^\frac12$; 
without loss of generality we restrict ourselves to $\vt$ such that $|a(\vt)|<1$. We can thus parametrize
the incidence relation as
\be\label{eqn zvt sim}
\zvt=\big\{ \left( x,\, x+r\om,\, x+t_1a(\vt)r\om + t_1b(\vt)r\ttom \right)\, : 
x\in\R^d,\, \om\in\sd,\, \ttom\in\sd\cap\, \om^\perp,\, r>0 \big\}.
\ee
Writing $\Psi=(\Psi_1,\Psi_2)$, one has
\bast
\wtnzt&=&\Big\{ \big(x, \tau_1 d_x\Psi_1+ \tau_2 d_x\Psi_2;\, x+r\om, \tau_1 d_y\Psi_1+ \tau_2 d_y\Psi_2;\\
& & \quad\qquad x+t_1a(\vt)r\om + t_1b(\vt)r\ttom, \tau_1 d_z\Psi_1+ \tau_2 d_z\Psi_2\big) \\
 & & \qquad\qquad :\, x\in\R^d,\, \om\in\sd,\, \ttom\in\sd\cap\, \om^\perp,\, r>0,\, \tau\in\R^2\setminus 0\, \Big\}.
\east
We  show that the projection $\pi_1:\wtnzt\to T^*\R^d$ onto the first factor  is a submersion, so that
Thm. \ref{thm kpoint old} applies with $\sigma=(1\, |\, 23)$ and no loss ($\beta=0$), yielding the threshold 
$\hd(E)>(1/3)(\max(d,2d)+2+0)=(2d+2)/3$, as claimed.
Strictly speaking we are not using Thm. \ref{thm kpoint new}, but will  in fact verify the nondegeneracy of $\pi_1$ 
separately on two open conic sets, $\{\tau_1\ne 0\}$ and $\{(t_1-a(\vt))\tau_1-t_2\tau_2\ne 0\}$, 
so that the microlocal approach is implicitly   being used.
\ms

Now,  one computes
\bast
d_x\Psi_1&=&|x-y|^{-1}|x-z|^{-1}(x-z)-|x-y|^{-3}|x-z|(x-y)\, ,\\
d_x\Psi_2&=&-|x-y|^{-3}|y-z|(x-y).
\east
Evaluating these at a point of $\zvt$ in terms of the variables in \eqref{eqn zvt sim}  using
\bast
\quad x-y= -r\om, & & x-z=-t_1a(\vt)r\om - t_1r\left(a(\vt)\om + b(\vt)\om'\right),  \\
|x-y|=r,\,\quad  & & |x-z|=t_1r, \quad |y-z|=t_2r,
\east
we find that

$$d_{x}(\Psi_1,\Psi_2)\big|_{\zvt}=\left(r^{-1}\left(\left(t_1-a(\vt)\right)\om -b(\vt)\ttom\right),\, 
-t_2r^{-1}\om \right).$$
Thus,
$$(x,\xi)=\pi_1(x,\om,\om',r,\tau)=\left(x, \tau_1r^{-1}\left(\left(t_1-a(\vt)\right)\om -b(\vt)\ttom\right)
-\tau_2t_2r^{-1}\om\right).$$

Since the spatial component of $\pi_1$ is the identity in $x$, to show that the canonical relation $C_{\vt}^{(1|23)}$ is 
nondegenerate, we need to show that $\rank D\xi/D(\om,\ttom,r,\tau)=d$ everywhere.
We have
$$D\xi/D\om = r^{-1}((t_1-a(\vt))\tau_1-t_2\tau_2)i_\om,\quad D\xi/D\tau_2=-t_2r^{-1}\om,$$
where $i_\om:T_\om\sd\hookrightarrow T_\om\R^d$ is the inclusion of the tangent space of the sphere at $\om$.
Thus, $D\xi/D(\om,\tau_2)$ is surjective on the set of $\tau$ such that $(t_1-a(\vt))\tau_1-t_2\tau_2\ne 0$.
\ms

On the other hand, 
$$D\xi/D\om'=-r^{-1}\tau_1b(\vt)j_{\om'},\qquad D\xi/D\tau_1
=r^{-1}\left(\left(t_1-a(\vt)\right)\om-b(\vt)\om'\right),$$
$$ \hbox{ and } \quad D\xi/D\tau_2=-t_2r^{-1}\om,$$
where $j_{\om'}:T_{\om'}\left(\sd\cap\, \om^\perp\right)\hookrightarrow T_{\om'}\R^d$ is the inclusion of the 
$(d-2)$-dimensional tangent space. As long as $\tau_1\ne 0$, this space and the two vectors that follow it span 
$\R^d$ (note that $b(\vt)\ne 0$)  and hence $D\xi/D(\om',\tau_1,\tau_2)$ is surjective. 
Since $\mathcal U_1=\{\tau_1\ne 0\}$ and $\mathcal U_2=\{(t_1-a(\vt))\tau_1-t_2\tau_2\ne 0\}$ form an open cover of 
$\wtnzt$, we are done.
\ms

As in the comment at the end of the previous section concerning congruences, 
it would be interesting to investigate whether the microlocal approach
can be applied to obtain results for similarities of $(k-1)$-simplices in $\R^d$.


\medskip

\end{document}